\pdfoutput=1

\documentclass[conference]{IEEEtran}
\usepackage{amssymb}
\usepackage{amsthm}
\theoremstyle{remark}
\usepackage{booktabs}
\usepackage{color}

\ifCLASSINFOpdf
 \usepackage[pdftex]{graphicx}
\graphicspath{{figures/}}
\else
\fi
\usepackage{amsmath}
\ifCLASSOPTIONcompsoc
  \usepackage[caption=false,font=normalsize,labelfont=sf,textfont=sf]{subfig}
\else
  \usepackage[caption=false,font=footnotesize]{subfig}
\fi

\makeatletter
\def\ps@IEEEtitlepagestyle{%
  \def\@oddfoot{\mycopyrightnotice}%
  \def\@evenfoot{}%
}
\def\mycopyrightnotice{%
  {\footnotesize Accepted to IEEE PES GM 2018. \copyright~IEEE 2018.\hfill}
  \gdef\mycopyrightnotice{}
}
\makeatother

\begin{document}
%
\title{The Value of Reactive Power \\for Voltage Control in Lossy Networks}

\author{
\IEEEauthorblockN{Matthew Deakin, \textit{Student Member, IEEE}, Thomas Morstyn, \textit{Member, IEEE},\\ Dimitra Apostolopoulou, \textit{Member, IEEE}, Malcolm McCulloch, \textit{Senior Member, IEEE}}
\IEEEauthorblockA{Department of Engineering Science, University of Oxford, Oxford, UK}
}


%


\maketitle

\begin{abstract}
Reactive power has been proposed as a method of voltage control for distribution networks, providing a means of increasing the amount of energy transferred from distributed generators to the bulk transmission network. The value of reactive power can therefore be measured according to an increase in \textit{transferred} energy, where the transferred energy is defined as the total generated energy, less the total network losses. If network losses are ignored, an error in the valuation of a given amount of reactive power will be observed (leading to reactive power provision being under- or over-valued). The non-linear analytic solution of a two-bus network is studied, and non-trivial upper and lower bounds are determined for this `valuation error'. The properties predicted by this two-bus network are demonstrated to hold on a three-phase unbalanced distribution test feeder with good accuracy. This allows for an analytic assessment of the importance of losses in the valuation of reactive power in arbitrary networks.
\end{abstract}


%
\IEEEpeerreviewmaketitle

\section{Introduction}

A drive to reduce reliance on fossil fuels has resulted in an unprecedented take-up of low-carbon technologies such as wind and solar. One of the features of this type of generation is that they are much more amenable to take-up within distribution networks. For example, in the UK, 54 \% of installed solar PV capacity is of size $< 5$ MW \cite{decc2017national} and so is unlikely to be connected directly to transmissions networks. On distribution networks, however, voltage rise tends to become a problem, eventually leading to curtailed real power \cite{masters2002voltage}.

Reactive power is often proposed as a means of controlling the voltage, and its tendency to increase network losses is well documented \cite{masters2002voltage}; however the impact of these losses across network types is not well understood (traditional loss analysis is focused on losses caused by load \cite{dickert2009energy}). For example, reactive power control was found to be the most cost effective network intervention in Germany \cite{stetz2013improved}, whilst the authors of \cite{gagrica2015microinverter} suggest that curtailment is the \textit{only} worthwhile method of voltage control in low voltage networks. In \cite{turitsyn2011options}, the authors discuss the natural trade-off between losses and voltage regulation but fail to recognise losses as a source of overestimation of the value of reactive power. Other studies that propose reactive power control fail to discuss losses at all \cite{sultan2015incorporating}, or, when undertaking economic analyses, do not equally weight the value of generated power and of losses equally \cite{stetz2014techno}. All of these studies fail to consider that an increase in real power generator export is only a valid utility function if the increase in losses caused by an intervention remains small, when compared to the increase in real power transferred to the grid (which is also potentially small). This motivates a search to investigate numerically the suitability of reactive power for voltage control and, particularly, how it varies across networks.

In this paper, the closed-form solution of the two bus power flow problem is studied to consider the real power losses that are caused by reactive power for voltage control. Thus, when ignored, these losses lead to a `valuation error': a difference between the increase in generated energy and the energy transferred to the bulk grid (given the installed reactive power). This leads to our main result: upper and (non-zero) lower bounds on this valuation error, which are independent of the temporal generation profile. Both under- and over- estimates of the value of reactive power are observed. The predicted properties are demonstrated on the unbalanced IEEE 34 bus distribution test feeder.

\section{Reactive Power for Voltage Control}

\begin{figure}
\centering
\includegraphics[width=0.32\textwidth]{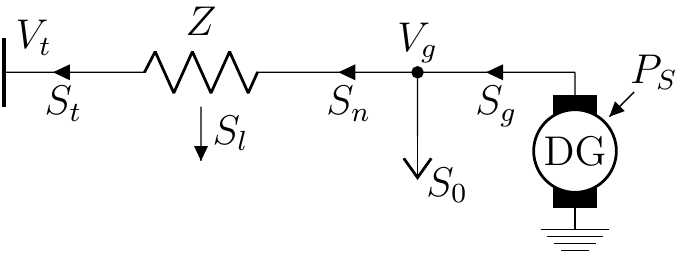}
\caption{Two bus power flow model (sign convention as indicated).}
\label{f:pesgm_trans}
\end{figure}

With reference to Fig. \ref{f:pesgm_trans}, we consider the impedance of the line $Z = R + \text{j}X$ where $R, X$ are both non-negative real numbers. We use $|\cdot |$ to denote the magnitude of a complex number. We use the notation $S_{(\cdot)} = P_{(\cdot)} + \text{j}Q_{(\cdot)}$ to represent apparent, real and reactive powers respectively. $S_{n},S_{t}$ represent the `net' and `transferred' apparent power, and $V_{g}$ the (complex) voltage at the generator bus. We choose $V_{t}$ as the reference bus, such that $V_{t}\in \mathbb{R}^{+}$. The generator is subject to a maximum upper voltage limit $V_{+} \in \mathbb{R}^{+}$.

The `net' power is given by $S_{n} = S_{g} - S_{0}$, and the `transferred' power by $S_{t} = S_{n} - S_{l}$, where $S_{0}$ represents a load, $S_{g}$ represents the power generated by distributed generator `DG', and $S_{l}$ represents the (complex) line losses. Finally, $P_{S}(\tau)$ represents the DG generation profile at time $\tau$.

\subsection{Two Bus Constant Voltage Load Flow}

It can be shown (e.g. from \cite{vournas2015maximum}) that, for a given $V_{t}$, that the lines of constant generator voltage $|V_{g}|$ lie on the circle
\begin{equation}\label{e:v2_circ}
\left( P_{n} - \dfrac{|V_{g}|^{2}R}{|Z|^{2}} \right) ^{2} + \left(  Q_{n} - \dfrac{|V_{g}|^{2}X}{|Z|^{2}} \right)^{2} = \left( \dfrac{V_{t}|V_{g}|}{|Z|} \right) ^{2} \,.
\end{equation}
For a given voltage limit we can therefore define the real and reactive power flows as
\begin{align}
P_{n}^{V} =& \dfrac{1}{|Z|^{2}}\left( R|V_{g}|^{2} \pm \sqrt{k_{ZV}W_{P}(Q_{n})} \right) \,, \label{e:PnQn}\\
Q_{n}^{V} =& \dfrac{1}{|Z|^{2}}\left( X|V_{g}|^{2} \pm \sqrt{k_{ZV}W_{Q}(P_{n})}\right) \,, \label{e:QnPn}
\end{align}
where
\begin{align}
W_{P}(Q_{n}) =& \left( \dfrac{R^{2}}{|Z|^{4}} -  \dfrac{Q_{n}^{2}}{|V_{g}|^{4}} - \dfrac{|V_{g}|^{2} - V_{0}^{2} - 2Q_{n}X}{|Z|^{2}|V_{g}|^{2}} \right)\,, \label{e:Wp} \\
W_{Q}(P_{n}) =& \left( \dfrac{X^{2}}{|Z|^{4}} -  \dfrac{P_{n}^{2}}{|V_{g}|^{4}} - \dfrac{|V_{g}|^{2} - V_{0}^{2} - 2P_{n}R}{|Z|^{2}|V_{g}|^{2}} \right)\,, \label{e:Wq} \\
k_{ZV} =& |Z|^{4}|V_{g}|^{4}\,.
\end{align}
By considering the quadrant of the circle we can determine the sign of $\pm$ in \eqref{e:PnQn} and \eqref{e:QnPn}. It is also useful to write down the identity \cite{deakin2017loss}
\begin{equation}
{P}_{l} = \dfrac{R}{|Z|^{2}} \left( V_{t}^{2} + 2(P_{n}R + Q_{n}X) - |V_{g}|^{2} \right) \,. \label{e:id_1}
\end{equation}
Prior to reaching the voltage limit, an increase in real power generated always yields an increase in transferred power. However, we re-state the result of \cite{deakin2017loss} that demonstrates that there exists a point on the voltage curve, the `marginal loss induced maximum power transfer point', where the losses increase at a greater rate than the increase in generated power. Therefore, beyond this point, real power should always be curtailed. This is therefore defined at the point
\begin{equation}\label{e:mlimpt_defn}
\dfrac{dP_{t}}{dP_{n}} = 0\,,
\end{equation}
and the reactive power at this point is given by
\begin{equation}\label{e:mlimpt}
Q_{n}^{\prime} = \dfrac{V_{+}^{2}X}{|Z|^{2}}\left( 1 - \dfrac{2V_{0}R}{V_{+}|Z|} \right) \,.
\end{equation}

\subsection{Operating Characteristic and the Value of Reactive Power}

\begin{figure}
\centering
\subfloat[]{\includegraphics[width=0.42\textwidth]{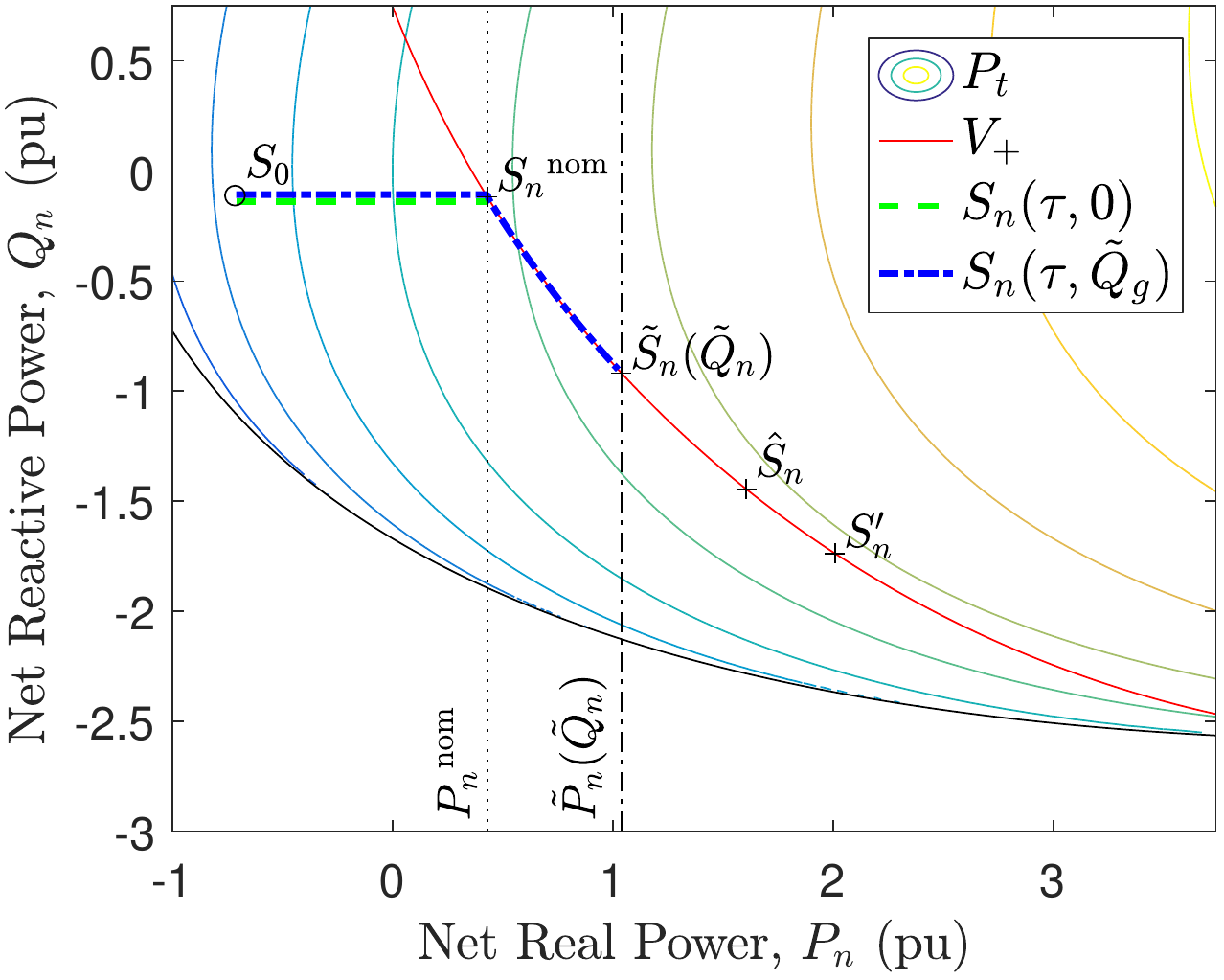}\label{f:analytic_method}}
\hfill
\subfloat[]{\includegraphics[width=0.42\textwidth]{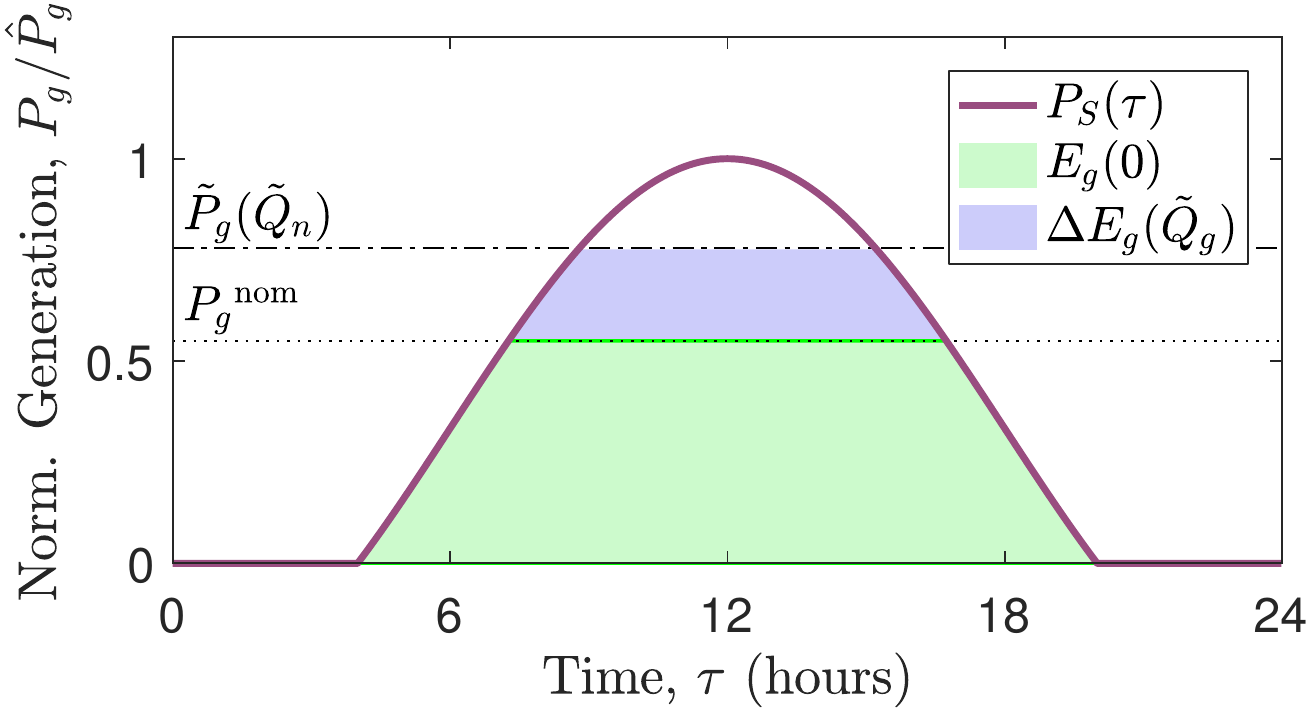}\label{f:solar_curves}}
\hfill
\caption{Operating characteristic \eqref{e:3stage_control} with respect to (a) the net power $P_{n}, Q_{n}$ through the impedance $Z$ and (b) the (normalised) generator power $P_{g}/\hat{P}_{g}$ in time. $S_{0}$ represents the feeder load, $S_{(\cdot)}^{\text{nom}}$ the nominal voltage-generation crossing, $\tilde{S}_{(\cdot)}$ the point at which generation is constrained (due to reactive power constraints), $\hat{S}_{(\cdot)}$ the generator peak output and $S_{(\cdot)}'$ the marginal loss-induced maximum power transfer point. Reactive power $\tilde{Q}_{g}$ increases the generated energy by $\Delta E_{g}$.}
\end{figure}

In this work a three-stage voltage control scheme is considered. Assume that exists a load $S_{0}$. As the generation $P_{S}$ is increased from zero, initially set $S_{g} = P_{S}$. However, as the real power is increased, the voltage curve will be hit at a power $P_{S} = P_{g}^{\text{nom}}$. Without access to reactive power, any power above this will be curtailed. However, if there exists a (limited) amount of reactive power available to the generator $\tilde{Q}_{g}$, we instead travel along the voltage circle defined by \eqref{e:v2_circ} until we reach this reactive power limit, at the real power $\tilde{P}_{g}$. Any addition real power will then be curtailed. This can be described as
\begin{equation}\label{e:3stage_control}
S_{g}(\tau,\tilde{Q}) = 
\begin{cases}
P_{S}(\tau) & \text{if } P_{S}(\tau) \leq P_{g}^{\text{nom}}\\
P_{S}(\tau) + \text{j}Q_{g}^{V}(P_{S}(\tau)) & \text{if } P_{g}^{\text{nom}} < P_{S}(\tau) \leq \tilde{P}_{g}\\
\tilde{P}_{g} + \text{j}\tilde{Q}_{g} & \text{if } P_{S}(\tau) > \tilde{P}_{g}\,.
\end{cases}
\end{equation}
using the identity $Q_{g}^{V} = Q_{n}^{V} + Q_{0}$. Fig. \ref{f:analytic_method} demonstrates this control in power space, while Fig. \ref{f:solar_curves} shows the operating characteristic as seen from the generator. 

The energy that can be generated $E_{g}$, and the energy that can be transferred $E_{t}$ (recalling $P_{t} = (P_{g} - P_{0}) - P_{l}$), are defined for a given reactive power limit $\tilde{Q}_{g}$ as
\begin{align*}
E_{g}(\tilde{Q}_{g}) = \int _{0}^{T} P_{g}(\tau,\tilde{Q}_{g}) \, d\tau \,, \\
E_{t}(\tilde{Q}_{g}) = \int _{0}^{T} P_{t}(\tau,\tilde{Q}_{g}) \, d\tau \,,
\end{align*}
where $T$ is the period over which the calculations are made. The reactive power results in an increase in energy generation and transfer as
\begin{align*}
\Delta E_{g}(\tilde{Q}_{g}) =& E_{g}(\tilde{Q}_{g}) - E_{g}(0)\,,\\
\Delta E_{t}(\tilde{Q}_{g}) =& E_{t}(\tilde{Q}_{g}) - E_{t}(0)\,.
\end{align*}

A na\"{i}ve approach to estimate the value of a given reactive power $\tilde{Q}_{g}$ might be to consider the increase in total generated energy, $\Delta E_{g}$ (the shaded area of Fig. \ref{f:solar_curves}). This is erroneous insofar as the increased flows result in a change in losses, which means that the value of the reactive power will be overestimated (or underestimated). Even if the power is generated at zero marginal cost, this still represents an error. The valuation `energy error' $\epsilon _{E}$ is therefore defined for a given generation profile $P_{S}(\tau)$ as
\begin{equation}\label{e:eE_defn1}
\epsilon _{E}(\tilde{Q}_{g}) = \dfrac{\Delta E_{g} - \Delta E_{t}}{\Delta E_{t}}\,.
\end{equation}
In general, this error will depend strongly on the generation profile $P_{S}(\tau)$.

We also define the instantaneous \textit{power} increase as
\begin{align*}
\Delta P_{g}(Q_{n}) =& P_{g}^{V}(Q_{n}) - P_{g}^{V}(0)\,,\\
\Delta P_{t}(Q_{n}) =& P_{t}^{V}(Q_{n}) - P_{t}^{V}(0)\,,
\end{align*}
and thus the valuation `power error' as
\begin{equation}\label{e:upe}
\epsilon _{P}(Q_{g}) = \dfrac{\Delta P_{g} - \Delta P_{t}}{\Delta P_{t}}\,.
\end{equation}
Note that, in contrast to $\epsilon_{E}$, $\epsilon_{P}$ is independent of the generation profile $P_{S}$, insofar as it is a defined for a given set of powers, rather that for a set of powers across time. 

Note that, with a slight abuse of notation, \eqref{e:eE_defn1} can be re-written as
\begin{equation}\label{e:uee_alternative}
\epsilon_{E}(\tilde{Q}_{g}) = \dfrac{\int _{0}^{T} \Delta P_{g} - \Delta P_{t} \, d\tau}{\int _{0}^{T} \Delta P_{t} \, d\tau}\,,
\end{equation}
where $\Delta P_{g},\Delta P_{t}$ would be determined at each point in time $\tau$, and additionally, assuming that we are on the voltage circle \eqref{e:v2_circ} for all $\tau$. In the case that $S_{n}$ does not change with $\tau$, we then observe that
\begin{equation}\label{e:eE_eP}
\epsilon_{E} = \epsilon_{P}\,.
\end{equation}

\section{Reactive Power Valuation Error Bounds}\label{s:utility_error_bounds}

In this section, we consider upper and lower bounds on the valuation power error $\epsilon_{P}$ as, by \eqref{e:eE_eP}, this represents the energy error that would be observed if we were to operate at either the maximum or minimum generation on the voltage curve \eqref{e:v2_circ} (i.e. at $P_{g}', P_{g}^{\text{nom}}$ respectively). In this section we assume no load such that $S_{g} = S_{n}$.

\subsection{Upper Power Error Bound}\label{s:upper_p_bnd}

\begin{figure}
\centering
\subfloat[]{\includegraphics[width=0.42\textwidth]{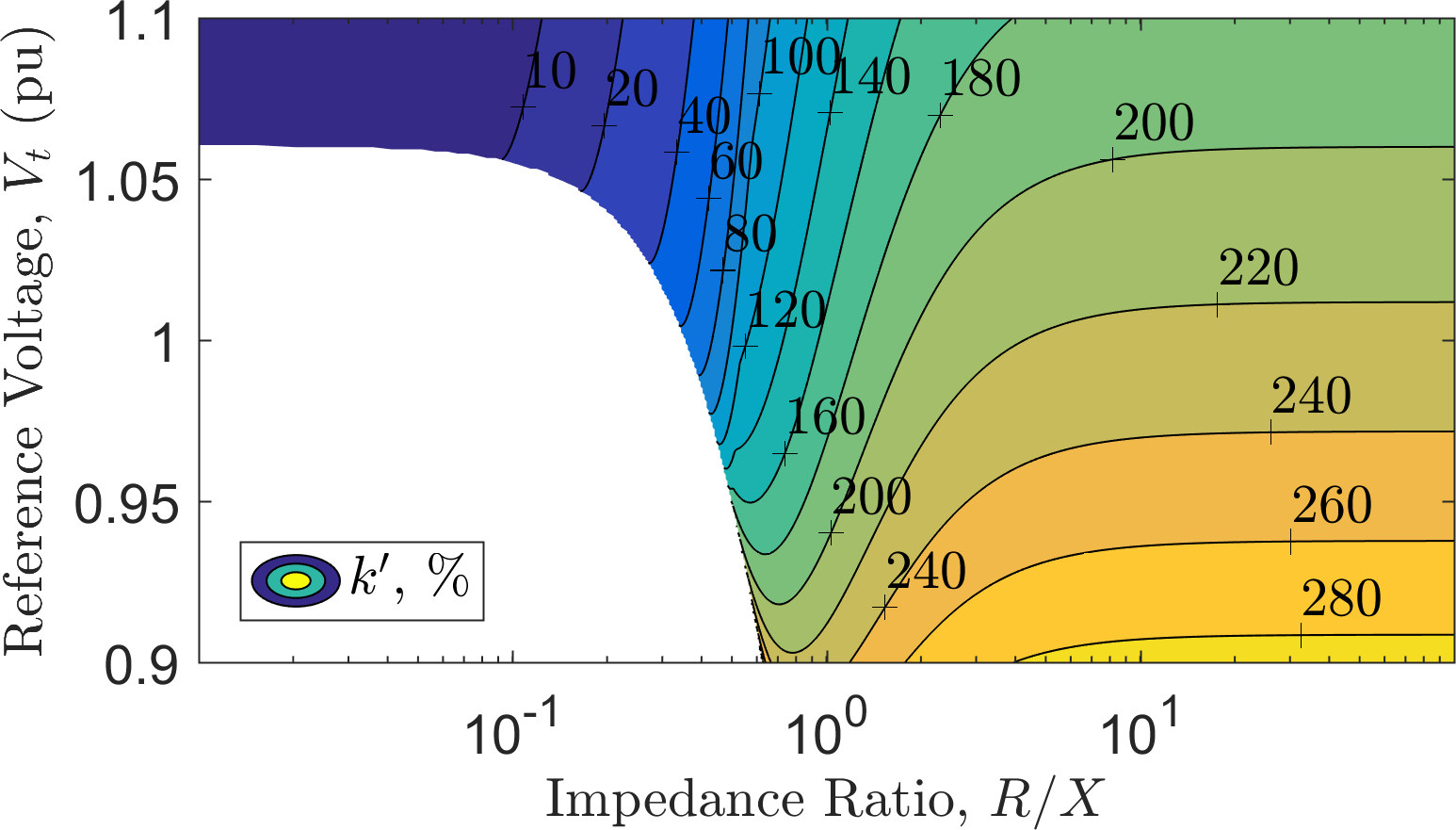}\label{f:k_dot}}
\hfill
\subfloat[]{\includegraphics[width=0.42\textwidth]{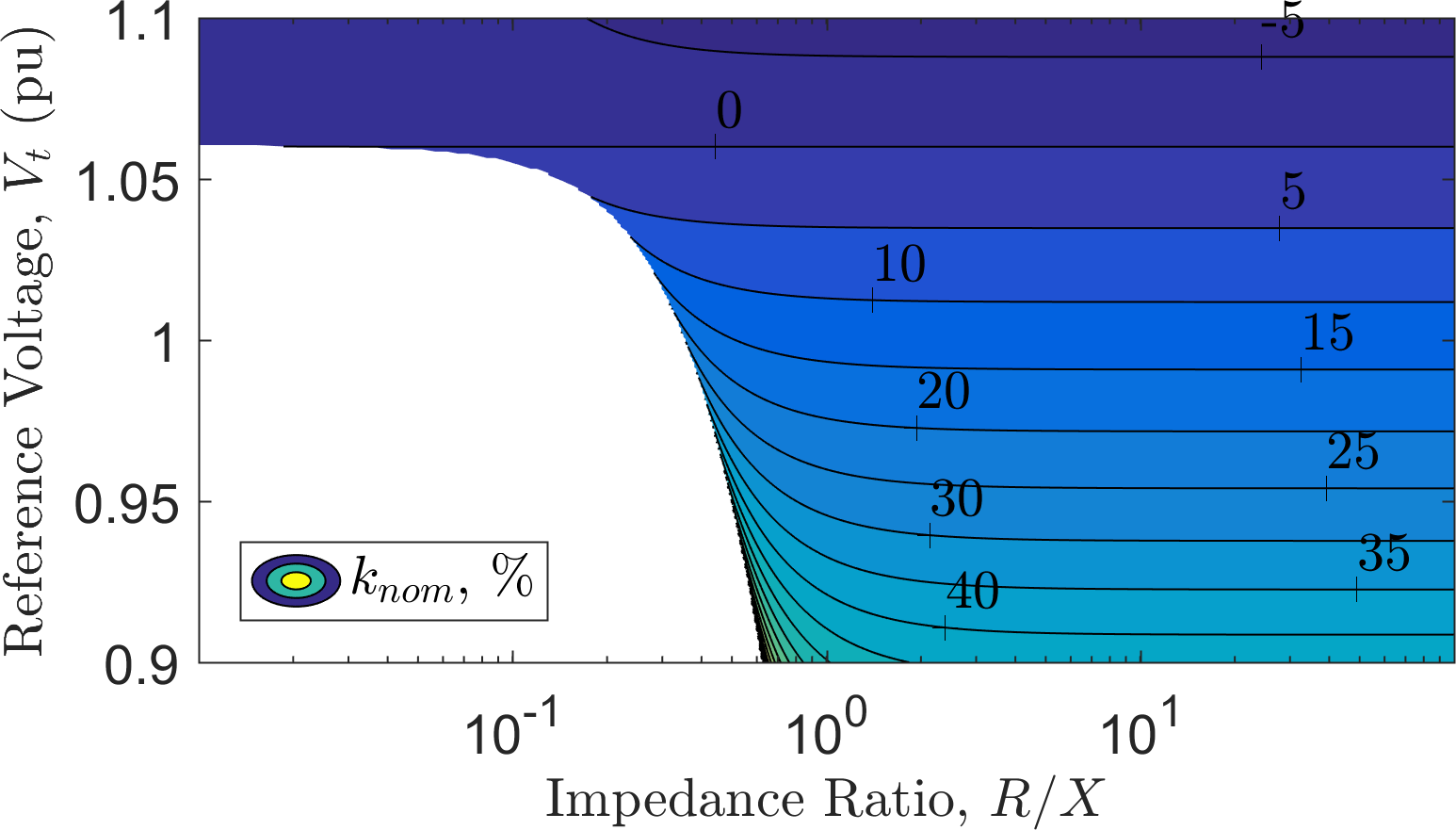}\label{f:k_bar}}
\hfill
\caption{Upper and lower bounds on the power error $\epsilon_{P}$, for $V_{+}=1.06$ pu: (a) $\max(\epsilon_{P}) = k'$, \eqref{e:kprm_defn}, and (b) $\min (\epsilon_{P}) = k_{nom}$, \eqref{e:knom_defn}. }
\end{figure}

To upper bound the power error $\epsilon_{P}$, we consider the case where the generator is operated at the marginal loss induced maximum power transfer point $P_{g}'$. We therefore define
\begin{equation}\label{e:kprm_defn}
k' = \epsilon_{P}(Q_{n}')\,,
\end{equation}
with $Q_{n}'$ defined as in \eqref{e:mlimpt}. The (closed form) solution to this is plotted in Fig. \ref{f:k_dot}. $\epsilon_{P}$ is undefined if the voltage curve does not intersect the line $Q_{n}=0$.

This bound represents the error that would be encountered if a generator was run at the marginal loss induced maximum power transfer point with 100\% availability. Even for conventional, dispatchable generators, this is not realistic. However, if information is known a priori about the maximum generator size $\hat{P}_{g}$ (due to, e.g., a thermal limit), then an alternative upper bound could be consider as $\hat{k} = \epsilon_{P}(\hat{P}_{n})$. Note that an operational constraint is always required to specify the maximum power that can be generated \cite{deakin2017loss}.

\subsection{Lower Power Error Bound}
First we note that the $\epsilon_{P}$ is undefined if $Q_{n} = 0$. Therefore, it is necessary to consider limiting cases to find a lower bound. As such, this lower bound, $k_{\text{nom}}$, is found as
\begin{equation}\label{e:knom_defn}
k_{\text{nom}} 	= \lim _{Q_{n} \to 0} \epsilon _{P}(Q_{n})\,.
\end{equation}
We rewrite \eqref{e:upe} as
\begin{equation*}
\epsilon _{P}(Q_{n}) = \dfrac{1}{\dfrac{\Delta P_{n}(Q_{n})}{P_{l}(Q_{n}) - P_{l}(0)} - 1}\,,
\end{equation*}
and define (using \eqref{e:PnQn})
\begin{align*}
k_{PQ} 	&= \lim _{Q_{n} \to 0} \dfrac{d\Delta P_{n}}{dQ_{n}}(Q_{n})\\
		&= \dfrac{-X|V_{g}|^{2}}{\sqrt{k_{ZV}W_{Q}(0)}}\,,
\end{align*}
Using these results we can thus re-write \eqref{e:knom_defn} as
\begin{align*}
k_{\text{nom}} &= \dfrac{1}{ k_{PQ}\lim_{Q_{n} \to 0} \left( \dfrac{Q_{n}}{P_{l}(Q_{n}) - P_{l}(0)}\right) - 1}\\
				&= \dfrac{dP_{l}}{dQ_{n}}(0)\dfrac{1}{k_{PQ} - \dfrac{dP_{l}}{dQ_{n}}(0)}\,,\\
\end{align*}
by the definition of the derivative of a function. Finally, differentiation of \eqref{e:id_1} yields
\begin{equation*}
k_{\text{nom}} = \dfrac{2R(Rk_{PQ} + X)}{2R(Rk_{PQ} + X) - |Z|^{2}k_{PQ}}\,.
\end{equation*}

This represents the error that is approached as we tend towards no reactive power. The value of this is plotted in Fig. \ref{f:k_bar}. Note that $\epsilon_{P}$ will \textit{not} generally tend to zero. To the contrary, we see that there are regions of the space where the limit tends to very large percentage errors (indeed, some values of $k_{\text{nom}}(Z,V_{t})$ are greater than some values of $k'(Z,V_{t})$). This implies that the losses will always represent a large percentage error, even if the magnitude of these losses are arbitrarily small.

Furthermore we see that, for some values of $(Z, V_{t})$ that $k_{\text{nom}} < 0$. This implies that ignoring losses could result in an underestimate of the value of reactive power. Given that $k' > 0$ in all cases, the magnitude of the marginal power error $|\epsilon_{P}|$ must therefore first decrease with $P_{n}$, before increasing.

\section{Case Study}

In this section we consider whether the properties predicted by analysis of the two bus case hold for a distribution feeder. Simulations are run in OpenDSS \cite{opendss2017} for the IEEE 34 bus distribution test feeder \cite{ieee2017distribution}. The code used to generate the results and figures of this paper is available at
\begin{center}
\texttt{https://github.com/deakinmt/pesgm18} .
\end{center}

\subsection{Network Analysis and Generation Profile}

\begin{table}[!t]
\caption{Network Parameters (pu). $S_{base} = 2.5$ MVA, $V_{base,LL} = 69$ kV.}
\label{t:cases}
\centering
\begin{tabular}{l l l l l l l}
\toprule
 $|Z|$ & $R/X$ & $V_{+}$ & $|S_{0}|$ & PF $(= P_{0}/|S_{0}|)$ \\
\midrule
0.203 & 1.85 & 1.06 	& 0.72 & 0.987 lagging \\
\bottomrule
\end{tabular}
\end{table}

We consider connecting a generator to bus 834 of the test feeder. This bus is situated towards the end of the feeder, as is a majority of the load. We calculate equivalent two bus network parameters as in \cite{deakin2017loss} (see Table \ref{t:cases}). The generator peak power is chosen as
\begin{equation}\label{e:generation_max}
\hat{P}_{g} = \left( P_{g}^{\text{nom}} + c( P_{g}' - P_{g}^{\text{nom}} ) \right) \,, \quad 0 < c \leq 1 \,.
\end{equation}
Over the course of a day ($T = 24$ hours) the generation profile is specified as
\begin{equation}\label{e:generation_profile}
P_{S}(\tau) = \max \left\{ \, \dfrac{\hat{P}_{g}}{3}\left(1 + 2\cos \left( \dfrac{\pi\tau}{12} \right) \right), 0\, \right\} \,,
\end{equation}
as in Fig. \ref{f:solar_curves} (which corresponds to a large horizontal array in July, with a clear sky, at a latitude of $48.8^{\circ}$, e.g., Paris \cite{kalogirou2009solar}). We utilise the control scheme described in \eqref{e:3stage_control} to determine the reactive powers.

\subsection{Valuation Power Error}\label{s:UPE}

To implement \eqref{e:3stage_control} in the feeder, a fine mesh  of real and reactive powers ($1000\times1200$) are specified at the generator at bus 834, with tap changers fixed for the simulations. For each point, the total network losses are calculated. Infeasible points are then removed (i.e. points where there is generator overvoltage). For a given $\tilde{Q}_{g}$, infeasible points are again removed (places where there is too much reactive power used). For each value of generated power $P_{g}$, the reactive power is then chosen to maximise the power transfer $P_{t}$. This has the effect of collapsing the mesh of real and reactive powers to a line, resembling the curves represented in Fig. \ref{f:analytic_method}.

\begin{figure}
\centering
\subfloat[]{\includegraphics[width=0.22\textwidth]{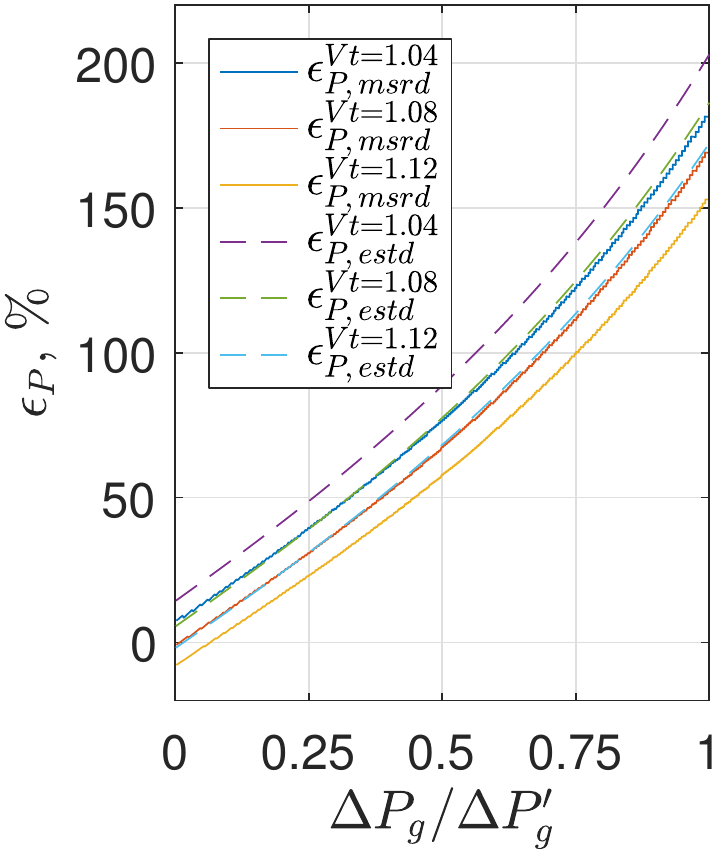}\label{f:upe}}
\hspace{1em}
\subfloat[]{\includegraphics[width=0.22\textwidth]{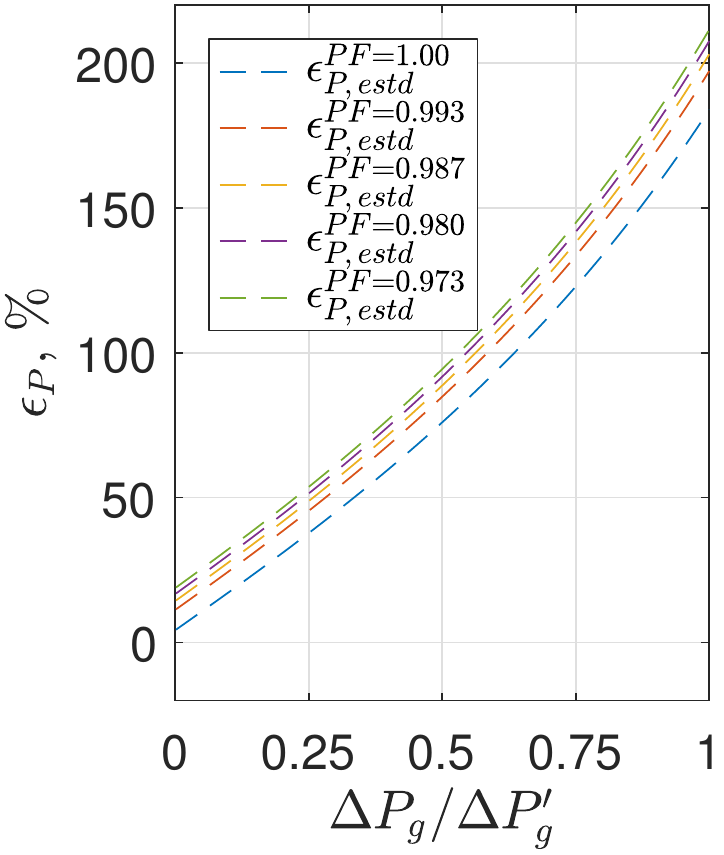}\label{f:nonlin_S0_pf}}
\hfill
\caption{The (a) measured and estimated valuation power error (b) estimated valuation power error sensitivity to (lagging) power factors ($V_{t}=1.04$).}
\end{figure}

These results are used to first derive the valuation power error $\epsilon_{P}$ for the network as a function of the increased real power $\Delta P_{g}$ (see Fig. \ref{f:upe}). We see that three properties predicted in Section \ref{s:utility_error_bounds} hold. In particular, we see (i) that there exists a non-zero lower bound for $\epsilon_{P}$, (ii) $\epsilon_{P}$ can sometimes become negative, and (iii) $\epsilon_{P}$ is upper bounded. The two bus approximation can be seen to broadly follow the error, although in this case it consistently over-estimates $\epsilon_{P}$. This over-estimate is presumed to be due to the voltage dependency of the reactive loads and capacitors in the full OpenDSS model. From Fig. \ref{f:nonlin_S0_pf}, we see $\epsilon_{P}$ is very sensitive to changes in load power factor.

\subsection{Valuation Energy Error}

To estimate the valuation energy error, $P_{n}^{\text{nom}}$ and $P'_{n}$ are first found for the network. A set of generation profiles $P_{S}(\tau)$ are created for three values of $c$, using \eqref{e:generation_max} and \eqref{e:generation_profile}. For each $\tau$, the real generator power closest to this value is found (using the curves from Section \ref{s:UPE}). This is then repeated for a range of values of the upper bound on reactive power $\tilde{Q}_{g}$. The analysis is then repeated for the two bus approximation, using the parameters of Table \ref{t:cases} (calculating separate two bus values of $P_{n}^{\text{nom}}$, $P'_{n}$, $P_{S}(\tau)$ etc).

\begin{figure}
\centering
\subfloat[]{\includegraphics[width=0.23\textwidth]{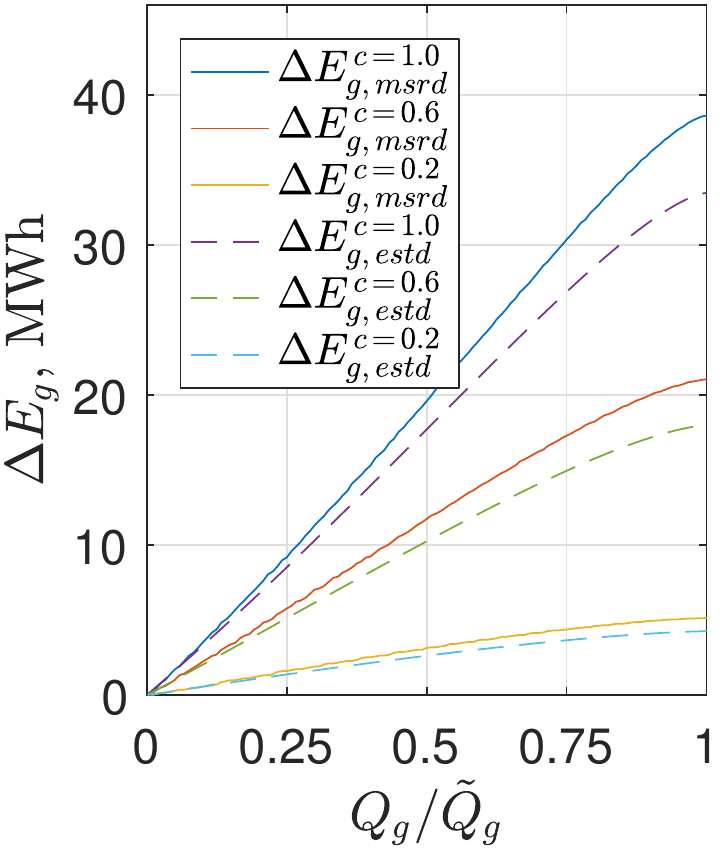}\label{f:eg_utility}}
\hspace{1em}
\subfloat[]{\includegraphics[width=0.23\textwidth]{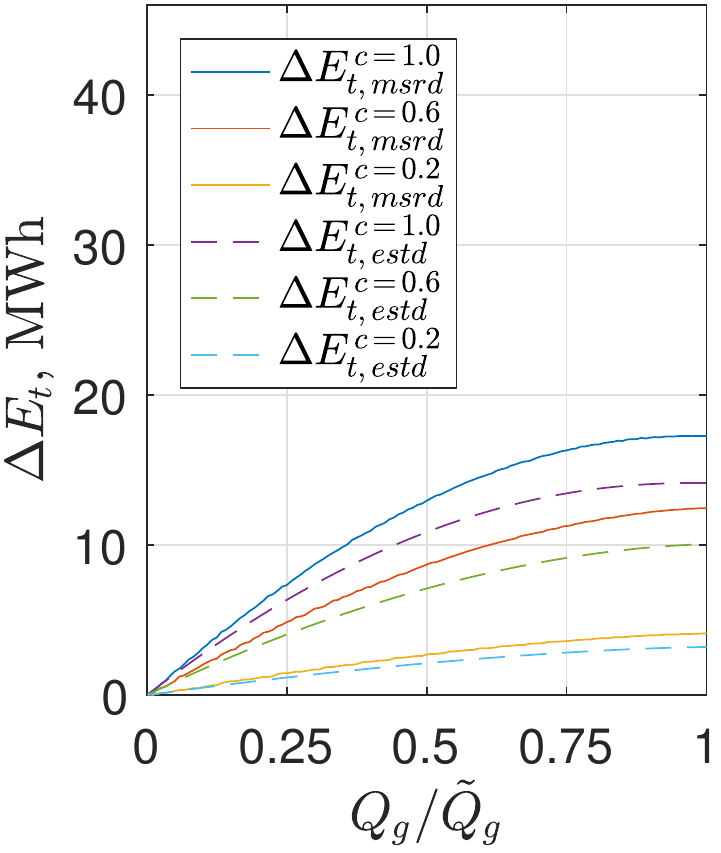}\label{f:et_utility}}
\hfill
\caption{The (a) generated and (b) transferred energy.}
\end{figure}

Figs. \ref{f:eg_utility} and \ref{f:et_utility} show the generator and transfer energy valuations $\Delta E_{g}$ and $\Delta E_{t}$ respectively, for $V_{t}=1.05$ pu. Generally, the measured and estimated results are of a similar magnitude and shape. A change in generation profile \eqref{e:generation_profile} changes the generated and transferred energy significantly, as expected. The two bus case this time underestimates these valuations. This is presumed to be due to the approximate nature of the two bus estimates of $P_{n}^{\text{nom}}, P_{n}'$ \cite{deakin2017loss}, and due to errors caused by the lumping of all of the load at the end of the feeder (as in Fig. \ref{f:pesgm_trans}).

\begin{figure}
\centering
\includegraphics[width=0.32\textwidth]{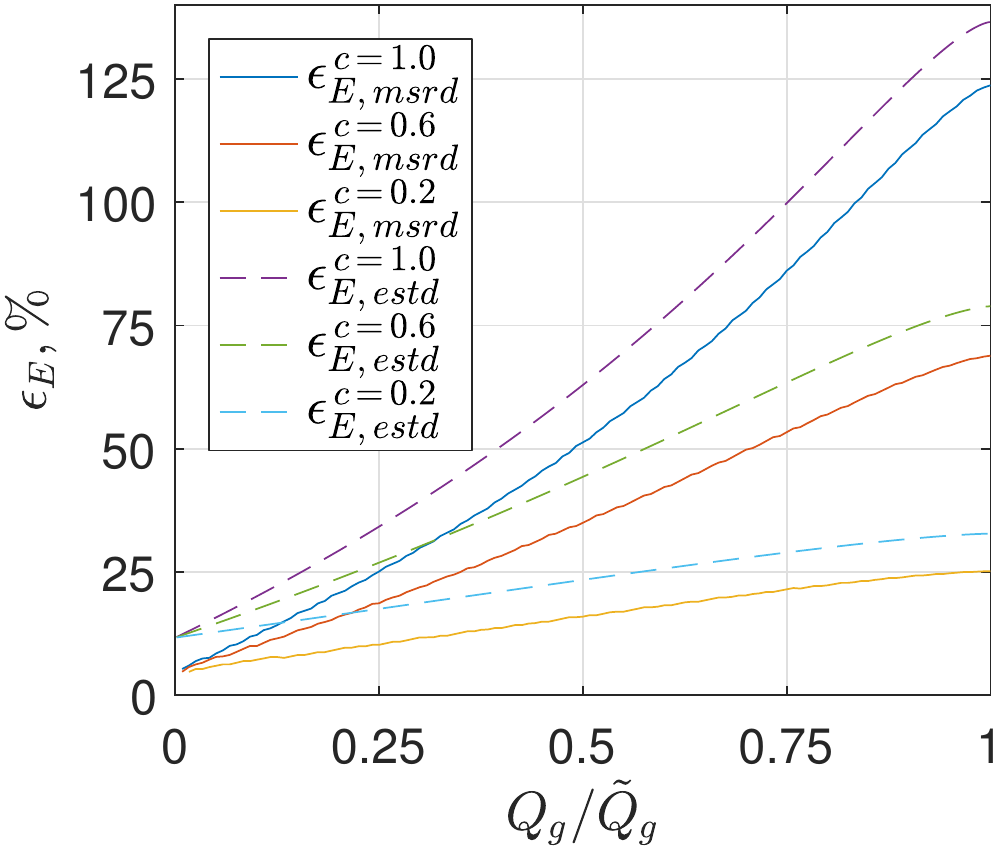}
\caption{The measured and estimated valuation energy error.}\label{f:uee}
\end{figure}

The valuation energy error $\epsilon_{E}$ is then plotted for these curves in Fig. \ref{f:uee}. It can be observed that there exists a non-zero lower bound on $\epsilon_{E}$. We see that this can be predicted by considering \eqref{e:uee_alternative} for the limiting case of an arbitrarily small amount of available reactive power. In addition, $\epsilon_{E}$ remains well below $k'$, as discussed in Section \ref{s:upper_p_bnd}. Here the two bus approximation overestimates the error, but the trend observed is accurate.

Finally, it is noted that, in practise, there are a wide range of criteria that must be satisfied when considering network interventions, including protection, thermal issues and feeder power factors. In addition, the load $S_{0}$ and grid voltage $V_{t}$ will change over time. Therefore, detailed simulations will always be necessary. However, we have shown that the \textit{properties} predicted by the analytic, computationally rapid two bus network appear to generally hold true, even though there is some error between the estimated and true energy transfer calculations.

\section{Conclusions}
The impact of losses on the measurement of the value of reactive power in distribution networks has been considered in this work. It has been shown analytically that, in general, ignoring losses leads to a non-zero error in calculations of the value of reactive power, even for arbitrarily small amounts of reactive power. An upper bound on the error, derived from a suitable maximum power transfer point, is given in closed form as a function of the network voltages and $R/X$ ratio. Finally, simulations using sample generation profiles on a test feeder demonstrate real networks exhibit the properties predicted by the two bus network. We conclude that the method presented can therefore be used to estimate the relative importance of losses in the determination of the value of reactive power, to aid network operators in the identification of the most promising locations for reactive power network interventions.

\section*{Acknowledgment}

The authors wish to thank the Oxford Martin Programme on Integrating Renewable Energy, the John Aird Scholarship, and the Clarendon Scholarship for their support.



%
\bibliographystyle{IEEEtran}
\bibliography{pesgm18_bib}{}

\begin{thebibliography}{10}
\providecommand{\url}[1]{#1}
\csname url@samestyle\endcsname
\providecommand{\newblock}{\relax}
\providecommand{\bibinfo}[2]{#2}
\providecommand{\BIBentrySTDinterwordspacing}{\spaceskip=0pt\relax}
\providecommand{\BIBentryALTinterwordstretchfactor}{4}
\providecommand{\BIBentryALTinterwordspacing}{\spaceskip=\fontdimen2\font plus
\BIBentryALTinterwordstretchfactor\fontdimen3\font minus
  \fontdimen4\font\relax}
\providecommand{\BIBforeignlanguage}[2]{{%
\expandafter\ifx\csname l@#1\endcsname\relax
\typeout{** WARNING: IEEEtran.bst: No hyphenation pattern has been}%
\typeout{** loaded for the language `#1'. Using the pattern for}%
\typeout{** the default language instead.}%
\else
\language=\csname l@#1\endcsname
\fi
#2}}
\providecommand{\BIBdecl}{\relax}
\BIBdecl

\bibitem{decc2017national}
{Department of Energy and Climate Change}, ``National statistics: Solar
  photovoltaics deployment in the uk (september 2017 update),''
  \url{https://www.gov.uk/government/statistics/solar-photovoltaics-deployment},
  accessed: 02-11-2017.

\bibitem{masters2002voltage}
C.~Masters, ``Voltage rise: the big issue when connecting embedded generation
  to long 11 kv overhead lines,'' \emph{Power engineering journal}, vol.~16,
  no.~1, pp. 5--12, 2002.

\bibitem{dickert2009energy}
J.~Dickert, M.~Hable, and P.~Schegner, ``Energy loss estimation in distribution
  networks for planning purposes,'' in \emph{PowerTech, 2009 IEEE Bucharest},
  2009.

\bibitem{stetz2013improved}
T.~Stetz, F.~Marten, and M.~Braun, ``Improved low voltage grid-integration of
  photovoltaic systems in germany,'' \emph{IEEE Transactions on sustainable
  energy}, vol.~4, no.~2, pp. 534--542, 2013.

\bibitem{gagrica2015microinverter}
O.~Gagrica, P.~H. Nguyen, W.~L. Kling, and T.~Uhl, ``Microinverter curtailment
  strategy for increasing photovoltaic penetration in low-voltage networks,''
  \emph{IEEE Transactions on Sustainable Energy}, vol.~6, no.~2, pp. 369--379,
  2015.

\bibitem{turitsyn2011options}
K.~Turitsyn, P.~{\v{S}}ulc, S.~Backhaus, and M.~Chertkov, ``{Options for
  control of reactive power by distributed photovoltaic generators},''
  \emph{Proceedings of the IEEE}, vol.~99, no.~6, pp. 1063--1073, 2011.

\bibitem{sultan2015incorporating}
S.~S. Sultan, V.~Khadkikar, and H.~H. Zeineldin, ``Incorporating pv inverter
  control schemes for planning active distribution networks,'' \emph{IEEE
  Transactions on Sustainable Energy}, vol.~6, no.~4, pp. 1224--1233, 2015.

\bibitem{stetz2014techno}
T.~Stetz, K.~Diwold, M.~Kraiczy, D.~Geibel, S.~Schmidt, and M.~Braun,
  ``Techno-economic assessment of voltage control strategies in low voltage
  grids,'' \emph{IEEE Transactions on Smart Grid}, vol.~5, no.~4, pp.
  2125--2132, 2014.

\bibitem{vournas2015maximum}
C.~Vournas, ``Maximum power transfer in the presence of network resistance,''
  \emph{IEEE Transactions on Power Systems}, vol.~30, no.~5, pp. 2826--2827,
  2015.

\bibitem{deakin2017loss}
M.~Deakin, T.~Morstyn, D.~Apostolopoulou, and M.~McCulloch, ``Loss induced
  maximum power transfer in distribution networks,'' Submitted to PSCC 2018.
  Preprint: \url{http://arxiv.org/abs/1710.07787}.

\bibitem{opendss2017}
EPRI, ``Opendss: Epri distribution system simulator,''
  \url{https://sourceforge.net/projects/electricdss/}, 2017.

\bibitem{ieee2017distribution}
{IEEE Power and Energy Society}, ``Distribution test feeders,''
  \url{https://www.ewh.ieee.org/soc/pes/dsacom/testfeeders/}, 2017.

\bibitem{kalogirou2009solar}
S.~A. Kalogirou, \emph{Solar Energy Engineering: Processes and Systems}.\hskip
  1em plus 0.5em minus 0.4em\relax Elsevier Science, 2009, ch.~2.

\end{thebibliography}

\end{document}